\documentclass[11pt,reqno]{amsart}
\usepackage[T1]{fontenc}
\usepackage[utf8]{inputenc}
\usepackage{lmodern}
\usepackage{amssymb,mathtools}
\usepackage{booktabs,microtype,enumitem,hyperref}
\hypersetup{hidelinks,pdftitle={Krzyz's conjecture for the seventh coefficient},pdfauthor={Xiaohui Xie}}
\setlist{itemsep=2pt,topsep=4pt}
\newtheorem{theorem}{Theorem}[section]
\newtheorem{lemma}[theorem]{Lemma}
\newtheorem{proposition}[theorem]{Proposition}
\theoremstyle{remark}

\theoremstyle{plain}
\newtheorem{corollary}[theorem]{Corollary}
\DeclareMathOperator{\Rea}{Re}
\DeclareMathOperator{\Ima}{Im}
\newcommand{\T}{\mathbb T}
\newcommand{\D}{\mathbb D}
\newcommand{\norm}[1]{\left\lVert#1\right\rVert}
\title[The seventh coefficient in the Krzy\.{z} conjecture]{Krzy\.{z}'s conjecture for the seventh coefficient}
\author{Xiaohui Xie}
\address{University of California, Irvine, Irvine, California, USA}
\email{xhx@uci.edu}
\date{October 30, 2026}

\begin{document}
\begin{abstract}
A holomorphic function bounded in modulus by one in the unit disk
has Taylor coefficients of modulus at most one. If the function
has no zeros, Krzy\.{z}'s conjecture predicts the stronger bound
$2/e$ for every coefficient of positive order. We establish this
bound for the seventh coefficient: $|a_7|\leq2/e$, with equality precisely
for the functions
$\eta\exp[-(1-\xi z^7)/(1+\xi z^7)]$, where $|\eta|=|\xi|=1$.
Thus the largest possible value of $|f^{(7)}(0)|$ in this class is
$7!\,2/e$, and the extremal functions are uniquely determined up to
rotations of the variable and multiplication by a constant of modulus one.
\end{abstract}
\maketitle

\section{Introduction}
Let $\D=\{z\in\mathbb C:|z|<1\}$ be the unit disk. A holomorphic
function on $\D$ has a Taylor expansion
\[
 f(z)=\sum_{k=0}^{\infty}a_kz^k,
 \qquad a_k=\frac{f^{(k)}(0)}{k!}.
\]
If $|f|\leq1$, Cauchy's integral formula gives $|a_n|\leq r^{-n}$
for every $0<r<1$. Letting $r$ tend to one yields $|a_n|\leq1$.
This bound is sharp: the function $f(z)=z^n$ has $a_n=1$.
It also has a zero. Krzy\.{z}'s problem asks how much smaller the
coefficient must be when zeros are forbidden throughout the disk.

Write
\[
 \mathcal B_*=
 \{f\text{ holomorphic on }\D:0<|f(z)|\leq1\text{ for }z\in\D\},
 \qquad
 M_n=\sup_{f\in\mathcal B_*}|a_n|.
\]
The conjecture is that $M_n=2/e$ for every $n\geq1$.
The lower bound $2/e$ holds at every index because the functions
\begin{equation}\label{eq:model}
 F_n(z)=\exp\!\left(-\frac{1-z^n}{1+z^n}\right)
       =e^{-1}+2e^{-1}z^n+O(z^{2n})
\end{equation}
all belong to $\mathcal B_*$. Indeed,
\[
 \Rea\frac{1-z^n}{1+z^n}
   =\frac{1-|z|^{2n}}{|1+z^n|^2}>0,
\]
so $|F_n(z)|<1$, while an exponential never vanishes.
Thus $M_n\geq2/e$. The conjecture asserts the reverse inequality
and predicts that the rotations of $F_n$ are the only extremal
functions. Here a rotation means replacing $f(z)$ by
$\eta f(\zeta z)$, where $|\eta|=|\zeta|=1$; this operation
preserves the moduli of all Taylor coefficients.

We establish the seventh-coefficient case.
\begin{theorem}\label{thm:seventh}
If $f(z)=\sum_{k\geq0}a_kz^k$ belongs to $\mathcal B_*$, then
$|a_7|\leq2/e$. Equality holds if and only if
\[
 f(z)=\eta\exp\!\left(-\frac{1-\xi z^7}{1+\xi z^7}\right),
 \qquad |\eta|=|\xi|=1.
\]
\end{theorem}

Krzy\.{z} posed the problem in 1968 \cite{Krzyz}. The first two
cases are classical. Hummel, Scheinberg, and Zalcman
\cite{HummelScheinbergZalcman} proved the third case; Tan \cite{Tan}
and Brown \cite{Brown} proved the fourth, and Samaris \cite{Samaris}
proved the fifth. Wikstr\"om \cite{Wikstrom} established the sixth
case. The survey of Mart\'in, Sawyer, Uriarte-Tuero, and Vukoti\'c
\cite{MartinEtAl} describes the extremal problem and many equivalent
formulations. MacDonald \cite{MacDonald} and Stupin \cite{Stupin}
give further necessary conditions for extremality. Related extremal
problems are studied by Instanes, Llinares, and Vukoti\'c
\cite{InstanesLlinaresVukotic}.

Agler and McCarthy \cite{AglerMcCarthy} related the coefficient problem to an
entropy inequality for polynomials whose zeros lie on the unit circle.
Lei and Zhang \cite{LeiZhang} proved that inequality and determined
its equality cases. Together, these results settle the coefficient
estimate for extremizers with the largest possible support. We
recall the representation that makes this distinction precise.

\subsection{The logarithm and a positive measure}
The condition that $f$ has no zeros allows us to take a holomorphic
logarithm on the simply connected disk. After multiplying a
nonconstant $f\in\mathcal B_*$ by a constant of modulus one, we may
assume that $f(0)>0$ and write
\[
 f=e^{-H},\qquad \Rea H=-\log|f|>0,
 \qquad H(0)=t=-\log f(0)>0.
\]
The strict inequality follows from the maximum principle.
We may therefore study the coefficients of $e^{-H}$ with $\Rea H>0$.

This change of variables already proves the first-coefficient bound.
The function
\[
 \omega(z)=\frac{H(z)-t}{H(z)+t}
\]
maps $\D$ into itself and satisfies $\omega(0)=0$. Schwarz's lemma
gives $|H'(0)|=2t|\omega'(0)|\leq2t$, and hence
\[
 |a_1|=|f'(0)|=e^{-t}|H'(0)|\leq2te^{-t}\leq2/e.
\]
For higher coefficients, several derivatives of $H$ enter the
expansion of $e^{-H}$. Their interaction is the main difficulty.

The Herglotz representation expresses a function with positive real
part as a positive superposition of elementary fractions:
\[
 H(z)=\int_{\T}\frac{1+\alpha z}{1-\alpha z}\,d\mu(\alpha),
 \qquad \mu(\T)=t,
 \qquad \T=\{\alpha\in\mathbb C:|\alpha|=1\}.
\]
Here $\mu$ is a finite positive measure. The real part of each
integrand is the Poisson kernel
$(1-|z|^2)/|1-\alpha z|^2$; thus the positivity of the measure
encodes the condition $\Rea H>0$.

An extremizer is a function for which $|a_n|=M_n$.
Section~\ref{sec:mass} proves that an extremizer exists and that its
representing measure is supported on at most $n$ points. Consequently,
after normalization, it has the form
\begin{equation}\label{eq:atoms}
\begin{aligned}
 f(z)&=\exp\!\left[-\sum_{j=1}^N\lambda_j
                         \frac{1+\alpha_jz}{1-\alpha_jz}\right],\\
 &\lambda_j>0,\qquad |\alpha_j|=1,\qquad
 t=\sum_{j=1}^N\lambda_j,\qquad N\leq n.
\end{aligned}
\end{equation}
The measure assigns mass $\lambda_j$ to the point $\alpha_j$;
such a point mass is called an atom. The support points of an
extremizer are distinct. The coefficient estimates below also allow
coincident points, which is useful when taking limits.
We call a support deficient when $N<n$.

The main step in the seventh-coefficient problem is the following
support theorem.
\begin{theorem}[Seven-atom property]\label{thm:fullsupport}
Every extremizer for $M_7$ has exactly seven atoms in its
Herglotz representation.
\end{theorem}
Once this is known, the entropy theorem gives the sharp coefficient
bound and all equality cases. Thus Theorem~\ref{thm:seventh}
follows from Theorem~\ref{thm:fullsupport} and the entropy reduction.
Wikstr\"om \cite{Wikstrom} used the same reduction for the sixth
coefficient, excluding smaller supports by bounds on masses and
support angles. Here the angle dependence is bounded uniformly,
so the remaining subdivisions involve only masses.

\subsection{Outline of the argument}
Section~\ref{sec:mass} establishes the finite-support representation
and the mass restrictions. Section~\ref{sec:entropy} then gives the
entropy identity and its equality consequence, explaining why it
remains to prove Theorem~\ref{thm:fullsupport}.
Every extremizer satisfies the mass restrictions
\begin{equation}\label{eq:intro-mass-bounds}
 \lambda_j\geq\frac1n,\qquad t\geq\frac nN.
\end{equation}
These restrictions hold for every coefficient index $n$.
For $n=7$, we prove the following strict estimate on the resulting
mass domains:
\begin{equation}\label{eq:strict-exclusion}
 N\leq6,\quad \lambda_j\geq1/7,\quad t\geq7/N
 \quad\Longrightarrow\quad |a_7|<\beta:=147/200.
\end{equation}
Since $\beta<2/e$, no such function can be an extremizer.
Indeed, bounding the tail of the exponential series by a geometric
series gives
$e<\sum_{k=0}^6 1/k!+1/(6\cdot6!)<87/32$, and
$147/200<64/87<2/e$.

To prove \eqref{eq:strict-exclusion}, we express $a_7$ as a polynomial
in the support parameters $\alpha_j$, multiplied by $e^{-t}$.
The Laurent-polynomial bounds in Section~\ref{sec:phase} hold for
all choices of these parameters at once. Section~\ref{sec:kernels}
provides bounds for the mass derivatives and a uniform estimate
when $t\geq69/20$. In the remaining bounded region,
Section~\ref{sec:cover} extends bounds at finitely many mass vectors
to the entire domain by interpolation. The precise finite
inequalities are collected in Proposition~\ref{prop:finite}.
Appendix~\ref{app:verification} describes the finite data and their
verification. The complete data and verification programs form
a separate Zenodo deposit.

The exclusion estimate proves Theorem~\ref{thm:fullsupport}.
Section~\ref{sec:conclusion} combines it with the entropy identity
to prove Theorem~\ref{thm:seventh}. For completeness,
Appendix~\ref{app:entropy} includes the simple-root case of the
Lei--Zhang inequality used here.

\subsection{Why it is enough to subdivide the masses}
An $N$-atom configuration has $N$ masses and $N$ angles.
A common rotation removes one angle without changing $|a_7|$,
leaving $2N-1$ real parameters. For six atoms this is an
eleven-dimensional domain.

Our point estimates bound the support parameters simultaneously.
At a fixed mass vector, a nonnegative Laurent polynomial gives an
upper bound valid for every $(\alpha_1,\ldots,\alpha_N)\in\T^N$.
The interpolation estimates preserve this uniformity as the
masses vary over a box. The finite cover therefore takes place
in at most six real dimensions.

The estimate \eqref{eq:strict-exclusion} applies to every configuration
satisfying the mass restrictions, whether or not it is stationary.
We therefore need no further stationarity equations in the subdivision.
The entropy argument treats the equality functions separately.

\subsection{Notation}
We write $dm=d\theta/(2\pi)$ for normalized arc measure on $\T$.
The symbol $[z^k]F$ denotes the coefficient of $z^k$ in $F$.
For an analytic polynomial $P$ and $p=1,2$, set
\[
 \|P\|_p=\left(\int_{\T}|P|^p\,dm\right)^{1/p}.
\]
Parseval's identity identifies $\|P\|_2$ with the Euclidean norm
of its coefficients, and Cauchy--Schwarz gives $\|P\|_1\leq\|P\|_2$.
For a holomorphic function, the corresponding Hardy norm is
$\|F\|_{H^p}=\sup_{0<r<1}\|F(r\,\cdot)\|_p$.
For polynomials these definitions agree.
Reflection relative to degree $d$ means
$P^*(z)=z^d\overline{P(1/\bar z)}$, with $d$ specified when it is
not the degree of $P$. A polynomial is self-inversive if
$P=cP^*$ for some $|c|=1$.
Every terminating decimal used in a bound denotes its exact
rational value.

\section{Extremizers and mass restrictions}\label{sec:mass}
We first reduce the extremal problem to finitely many positive masses.
The finite-support representation is classical
\cite{HummelScheinbergZalcman,MartinEtAl}; the derivation below also
fixes the notation for the mass estimates. Throughout this section,
$n\geq1$ is arbitrary. Multiplication by a constant of modulus one
and rotation of the variable allow us to assume that $a_0>0$ and
$a_n>0$.

\subsection{Existence and finite support}
Let $(f_m)$ be a maximizing sequence. Montel's theorem gives a
subsequence converging uniformly on compact subsets of $\D$ to a
holomorphic function $f$. Taylor coefficients converge under this
limit, so $|a_n(f)|=M_n\geq2/e>0$. In particular, $f$ is not
identically zero. By Hurwitz's theorem, it has no zeros, and the
maximum principle gives $|f|<1$. Thus $f$ is a global extremizer.

For a normalized extremizer, write $f=e^{-H}$ with $H(0)=t>0$.
Its Herglotz representation is
\[
 H(z)=\int_{\T}\frac{1+\alpha z}{1-\alpha z}\,d\mu(\alpha),
 \qquad \mu(\T)=t.
\]
Define
\begin{equation}\label{eq:P}
 P(w)=a_n+2\sum_{k=1}^n a_{n-k}w^k,
 \qquad h(\psi)=\Rea P(e^{i\psi}).
\end{equation}
Adding mass $s\geq0$ at $w\in\T$ replaces $f$ by
$f(z)\exp[-s(1+wz)/(1-wz)]$. The derivative of its $n$th
coefficient at $s=0$ is $-P(w)$. Since the real part of that
coefficient cannot increase from $a_n=M_n$, we obtain $h\geq0$.
On the other hand, replacing $\mu$ by $(1+s)\mu$ is admissible
for both signs of sufficiently small $s$. Stationarity therefore
gives $\int h\,d\mu=0$. Positivity now implies that $\mu$ is
supported on the zeros of $h$.

The trigonometric polynomial $h$ has degree exactly $n$, because
its top Fourier coefficient is $a_0>0$. A nonnegative trigonometric
polynomial of degree $n$ has at most $n$ distinct zeros on the
circle: each zero has even multiplicity, and multiplication by
$e^{in\psi}$ gives a polynomial of degree $2n$.
This proves \eqref{eq:atoms} with $1\leq N\leq n$ and distinct
$\alpha_j$. The variations above concern a fixed Taylor coefficient;
they may all be justified by Cauchy's formula on a circle of radius
less than one.

By the Fej\'er--Riesz theorem, there is a polynomial $q$ of degree
$n$ such that
\[
 h(\psi)=|q(e^{i\psi})|^2.
\]
Thus every support point $\alpha_j$ is a zero of $q$. We next show
that these are precisely its unit-circle zeros and that they are
simple.

\subsection{A stationary phase identity}
Put
\[
 \Phi(\psi)=\sum_{j=1}^N
       \frac{2\lambda_j}{|e^{i\psi}-\alpha_j|^2}.
\]
Away from the support points, one has
\begin{equation}\label{eq:phase}
 \Phi(\psi)h(\psi)
   =nh(\psi)-\Rea[e^{i\psi}P'(e^{i\psi})].
\end{equation}
This is the phase identity of
\cite[Theorem~2 and Corollary~2.1]{MacDonald} and
\cite[Lemma~13.2]{Stupin}. We give a derivation.

Let $g=\log f$ and define
\[
 R(z)=\sum_{k=0}^n a_kz^k
           +\sum_{k=0}^{n-1}\overline{a_k}z^{2n-k}.
\]
For $z=e^{i\theta}$,
$R(z)=z^nh(-\theta)$. Its zeros at $\overline{\alpha_j}$ have
order at least two. They cancel the double poles of $zg'(z)$,
so $S(z)=zg'(z)R(z)$ is a polynomial. Moreover,
$zg'(z)=O(1/z)$ at infinity, giving $\deg S\leq2n-1$.
On the circle,
\[
 zg'(z)=\Phi(-\theta)\in\mathbb R
\]
away from the poles. Hence $R=R^*$ and $S=S^*$ when reflection
is taken relative to degree $2n$.

Since $R=f+O(z^{n+1})$, the coefficients of $S$ through degree
$n$ are $ka_k$. Reflection determines the remaining coefficients.
It follows that
\begin{align*}
 e^{-in\theta}S(e^{i\theta})
 &=na_n+2\Rea\sum_{k=0}^{n-1}ka_ke^{i(k-n)\theta}\\
 &=nh(-\theta)-\Rea[e^{-i\theta}P'(e^{-i\theta})].
\end{align*}
The left side also equals $\Phi(-\theta)h(-\theta)$, which
proves \eqref{eq:phase}.

\subsection{The masses at the support points}
We need a boundary derivative formula. Let $D_h(\zeta)$ denote
the outward radial derivative of the harmonic extension of $h$ at
$\zeta\in\T$. Thus
\[
 D_h(\zeta)=\Rea(\zeta P'(\zeta)).
\]
This is the derivative of the Poisson extension of the boundary
function $|q|^2$; it is not the radial derivative of $|q(z)|^2$
inside the disk.

\begin{lemma}\label{lem:normal-derivative}
Write $q(z)=Q+(z-1)r(z)$, where $Q=q(1)$. Then
\begin{equation}\label{eq:normal-derivative}
 D_h(1)=2\Rea(\overline Q r(1))-2\|r\|_2^2.
\end{equation}
In particular, if $q(\alpha)=0$ and
$r_\alpha(z)=q(z)/(z-\alpha)$, then
$D_h(\alpha)=-2\|r_\alpha\|_2^2<0$.
\end{lemma}
\begin{proof}
The harmonic extension of $e^{ik\theta}$ is
$r^{|k|}e^{ik\theta}$, so the radial derivative has Fourier
multiplier $|k|$. On the unit circle,
\[
 |q|^2=|Q|^2+2\Rea\bigl(\overline Q(z-1)r(z)\bigr)
                    +(2-z-z^{-1})|r(z)|^2.
\]
The second term contributes $2\Rea(\overline Q r(1))$ at $z=1$.
For the last term, the identity
\[
 2|k|-|k+1|-|k-1|=-2\mathbf1_{\{k=0\}}
\]
leaves only the constant Fourier coefficient of $|r|^2$,
which is $\|r\|_2^2$. This proves \eqref{eq:normal-derivative}.
The formula at $\alpha$ follows by rotating $\alpha$ to one.
\end{proof}

\begin{lemma}[Individual masses]\label{lem:individual}
The unit-circle zeros of $q$ are exactly the points $\alpha_j$,
and each is simple. If $r_j(w)=q(w)/(w-\alpha_j)$, then
\begin{equation}\label{eq:mass}
 \lambda_j=\frac{\|r_j\|_2^2}{|r_j(\alpha_j)|^2}
            \geq\frac1n.
\end{equation}
\end{lemma}
\begin{proof}
At a unit-circle zero $\alpha$, Lemma~\ref{lem:normal-derivative}
gives $-D_h(\alpha)>0$. If $\alpha$ were not a support point,
the left side of \eqref{eq:phase} would tend to zero at $\alpha$,
whereas the right side would tend to $-D_h(\alpha)>0$.
If a support zero had multiplicity at least two in $q$, then $h$
would vanish to order at least four. Its product with the double
pole of $\Phi$ would again tend to zero, giving the same contradiction.
Thus all unit-circle zeros are simple and belong to the support.
This conclusion is also obtained in \cite[Theorem~13.2]{Stupin}.

For $\psi=\arg\alpha_j+\delta$, we have
\[
 h(\psi)=|r_j(\alpha_j)|^2\delta^2+O(\delta^3),
 \qquad
 \Phi(\psi)=2\lambda_j\delta^{-2}+O(1).
\]
Taking limits in \eqref{eq:phase} proves the identity in
\eqref{eq:mass}. The polynomial $r_j$ has degree at most $n-1$,
so Cauchy--Schwarz applied to its $n$ coefficients gives
$|r_j(\alpha_j)|^2\leq n\|r_j\|_2^2$.
\end{proof}

\subsection{The total mass}
The individual bounds give $t\geq N/n$. The stationary identity
also gives the complementary bound needed for small supports.
\begin{lemma}[Total mass]\label{lem:total}
A normalized $n$th-coefficient extremizer with $N$ support points
satisfies $t\geq n/N$.
\end{lemma}
\begin{proof}
Choose a point $\zeta\in\T$ at which $q$ does not vanish. Rotate
the boundary coordinate so that $\zeta$ becomes one, carrying $h$,
$q$, and the support points along with it. The phase identity is
unchanged by this rotation. Write
\[
 q(z)=Q+(z-1)r(z),\qquad
 r(z)=\sum_{\ell=0}^{n-1}r_\ell z^\ell,
 \qquad Q\ne0.
\]
Substituting \eqref{eq:normal-derivative} into \eqref{eq:phase}
and completing squares gives
\begin{equation}\label{eq:complete}
 \Phi(0)=\frac n2+
   2\sum_{\ell=0}^{n-1}\left|\frac{r_\ell}{Q}-\frac12\right|^2
   \geq\frac n2.
\end{equation}
Since the evaluation point was arbitrary, $\Phi\geq n/2$
everywhere away from its poles.

Now consider the rational function $H=-\log f$. Its poles are
$\overline{\alpha_j}$. On each open arc between consecutive poles,
$H(e^{i\theta})$ is purely imaginary and
\[
 \frac{d}{d\theta}\Ima H(e^{i\theta})=-\Phi(-\theta)<0.
\]
The imaginary part decreases from $+\infty$ to $-\infty$ on each
arc. Hence $H$ has one simple zero $\beta_\ell$ on each arc,
for $1\leq\ell\leq N$. Its numerator has degree at most $N$,
so these are all its zeros.

The residue of $1/H$ at $\beta_\ell$ is
$-\beta_\ell/\Phi(-\arg\beta_\ell)$, because
$\beta_\ell H'(\beta_\ell)=-\Phi(-\arg\beta_\ell)$.
Thus the corresponding partial-fraction mass is
\[
 \nu_\ell=\frac{1}{2\Phi(-\arg\beta_\ell)}\leq\frac1n.
\]
The difference
\[
 \frac1H-\sum_{\ell=1}^N
       \nu_\ell\frac{\beta_\ell+z}{\beta_\ell-z}
\]
has no finite poles and is finite at infinity, since
$H(\infty)=-t\ne0$. It is therefore constant. Its boundary values
away from the poles are purely imaginary, while its value at zero
is real; the constant must be zero. Evaluating at zero yields
\[
 \frac1t=\sum_{\ell=1}^N\nu_\ell\leq\frac Nn,
\]
as required.
\end{proof}

\section{Full support and the equality case}\label{sec:entropy}
For a polynomial $p$ with $\|p\|_2=1$, the function $|p|^2$ has
integral one on $\T$. Its entropy functional is
\[
 \mathcal E(p)=\int_{\T}|p|^2\log|p|^2\,dm,
 \qquad 0\log0=0.
\]
Agler and McCarthy \cite{AglerMcCarthy} showed that this functional
determines the extremal coefficient when the representing measure
has exactly $n$ atoms. We prove the identity in our notation.
\begin{lemma}[Full-degree entropy identity]\label{lem:entropyidentity}
Let a normalized global $n$th-coefficient extremizer have exactly
$n$ distinct positive atoms. Put
$B(z)=\prod_{j=1}^n(1-\alpha_jz)=\sum_{j=0}^n b_jz^j$ and
$A=\sum_{j=0}^n|b_j|^2$. Then
\[
 a_n=e^{-t}A,\qquad
 \mathcal E(B/\sqrt A)=t-\log A=-\log a_n.
\]
\end{lemma}
\begin{proof}
The polynomial $R$ defined in the phase-identity derivation has
degree $2n$, constant coefficient $a_0=e^{-t}$, and a double zero
at each $\overline{\alpha_j}$. Hence $R=e^{-t}B^2$.
Since all zeros of $B$ lie on $\T$ and $B(0)=1$, its coefficients
satisfy $b_{n-j}=b_n\overline{b_j}$ and $|b_n|=1$.
Consequently $[z^n]B^2=b_nA$.
Comparing the middle coefficient of $R=e^{-t}B^2$, and using
$a_n>0$, gives $b_n=1$ and $a_n=e^{-t}A$.
Define $m_k=\sum_j\lambda_j\alpha_j^k$ and
$s_k=\sum_j\alpha_j^k$. Since
\[
 \log(f/e^{-t})=-2\sum_{k\geq1}m_kz^k,\qquad
 2\log B=-2\sum_{k\geq1}(s_k/k)z^k,
\]
comparison of logarithms through degree $n$ in
$f=e^{-t}B^2+O(z^{n+1})$ gives $m_k=s_k/k$ for $1\leq k\leq n$.

Let $C_k=\sum_{r=0}^{n-k}b_{r+k}\overline{b_r}$. Then
\[
 |B(e^{i\theta})|^2=A+2\Rea\sum_{k=1}^n C_ke^{ik\theta}.
\]
The measure assigning mass $\lambda_j$ to $\overline{\alpha_j}$
is supported on the zeros of $B$. Integrating the preceding
identity against this measure gives
\[
 tA+2\Rea\sum_{k=1}^n m_k\overline{C_k}=0.
\]
The Fourier expansion of $\log|B|^2$, paired with the finite
polynomial $|B|^2$, gives
\[
 \int_{\T}|B|^2\log|B|^2\,dm
 =-2\Rea\sum_{k=1}^n\frac{s_k}{k}\overline{C_k}=tA.
\]
This pairing is justified by radial approximation: each unit-root
logarithm converges in $L^1$, and $|B|^2$ is bounded.
Dividing by $A$ and subtracting $\log A$ proves the result.
\end{proof}
The simple-root case of the Lei--Zhang entropy theorem is
\begin{equation}\label{eq:entropybound}
 \mathcal E(p)\geq1-\log2,
 \qquad \|p\|_2=1,
\end{equation}
for every polynomial of degree $n\geq2$ with distinct zeros on $\T$.
Equality holds precisely for
$p(z)=c(\omega+z^n)$, $|\omega|=1$, $|c|=1/\sqrt2$.
Appendix~\ref{app:entropy} reproduces their argument in this case.
Lemma~\ref{lem:individual} ensures that the extremal polynomial
$B$ satisfies this simplicity assumption.

\begin{corollary}[Full-support reduction]\label{cor:fullsupport}
Let $n\geq2$. If an extremizer for $M_n$ has exactly $n$ atoms,
then $M_n=2/e$, and this extremizer is a rotation of $F_n$.
\end{corollary}
\begin{proof}
Normalize the extremizer so that $a_0>0$ and $a_n>0$.
Lemma~\ref{lem:entropyidentity} and \eqref{eq:entropybound} give
$-\log a_n\geq1-\log2$, hence $a_n\leq2/e$. Since
$a_n=M_n\geq2/e$, equality holds.

The equality statement in the entropy inequality forces all
intermediate coefficients of $B$ to vanish. Since $B(0)=1$ and
$b_n=1$, we have $B=1+z^n$.
The logarithmic moments in the proof of
Lemma~\ref{lem:entropyidentity} therefore satisfy $m_k=0$ for
$1\leq k<n$ and $m_n=-1$. The support parameters are the $n$th
roots of $-1$. The discrete Fourier transform on these points
gives equal masses $t/n$, and $m_n=-t$ gives $t=1$. Thus the
normalized extremizer is $F_n$.
Undoing the input and output rotations proves the assertion.
\end{proof}

\section{Bounds uniform in the support parameters}\label{sec:phase}
For fixed masses, the seventh coefficient is a finite polynomial in
$\alpha_1,\ldots,\alpha_N$. To write it explicitly, define $p_0(x)=1$
and
\begin{equation}\label{eq:pk}
 p_k(x)=\sum_{r=1}^k\binom{k-1}{r-1}\frac{(-2x)^r}{r!},
 \qquad k\geq1.
\end{equation}
Expanding the exponential and then $(1-z)^{-r}$ shows that
\[
 \exp\!\left(-\frac{2xz}{1-z}\right)=\sum_{k=0}^{\infty}p_k(x)z^k.
\]
For a multi-index $k=(k_1,\ldots,k_N)$ of nonnegative integers,
write $|k|=k_1+\cdots+k_N$ and
$\alpha^k=\alpha_1^{k_1}\cdots\alpha_N^{k_N}$. Multiplication of
the $N$ power series in \eqref{eq:atoms} gives
\begin{equation}\label{eq:phasepoly}
 a_7=e^{-t}F_\lambda(\alpha),\qquad
 F_\lambda(\alpha)=\sum_{|k|=7}
       \left(\prod_{j=1}^Np_{k_j}(\lambda_j)\right)\alpha^k.
\end{equation}
All its coefficients are rational when the masses are rational.

A common rotation of the support parameters multiplies
$F_\lambda$ by $e^{7i\theta}$. Therefore
\begin{equation}\label{eq:phase-rotation}
 \sup_{\T^N}|F_\lambda|=\sup_{\T^N}\Rea F_\lambda.
\end{equation}
It is thus enough to majorize the real part uniformly on the torus.
We do so by adding a nonnegative Laurent polynomial.

Consider the monomial vector
\[
 v(\alpha)=\bigl((\alpha^{-k})_{|k|=3},
                         (\alpha^\ell)_{|\ell|=4}\bigr).
\]
Its length is $\binom{N+2}{3}+\binom{N+3}{4}$.
The cross terms between these two groups contain the monomials of
total degree seven occurring in $F_\lambda$.
For an integer matrix $L$ with this many rows, an integer $\kappa$, and
a positive rational number $s$, put
\[
 Q(\alpha)=s2^{2\kappa}\|L^Tv(\alpha)\|_2^2\geq0.
\]

\begin{lemma}[A Laurent-polynomial majorant]\label{lem:phasebound}
Write
\[
 \Rea F_\lambda+Q=\sum_{\nu\in\mathbb Z^N}R_\nu\alpha^\nu,
 \qquad
 U=R_0+\sum_{\nu\ne0}|R_\nu|.
\]
Then
\begin{equation}\label{eq:phasebound}
 \sup_{\alpha\in\T^N}|a_7|\leq e^{-t}U.
\end{equation}
In particular, if $b>0$ and
\begin{equation}\label{eq:pointaccept}
 U<bE_{128}(t),\qquad
 E_{128}(t)=\sum_{k=0}^{128}\frac{t^k}{k!},
\end{equation}
then $|a_7|<b$ for every choice of support parameters.
\end{lemma}
\begin{proof}
On $\T^N$, positivity of $Q$ and the triangle inequality give
\[
 \Rea F_\lambda\leq\Rea F_\lambda+Q
                  \leq R_0+\sum_{\nu\ne0}|R_\nu|=U.
\]
Equation~\eqref{eq:phase-rotation} proves \eqref{eq:phasebound}.
Finally, $E_{128}(t)<e^t$ for $t>0$, so
$e^{-t}U<be^{-t}E_{128}(t)<b$.
\end{proof}

The sum defining $U$ includes every Fourier coefficient, including
those introduced by $Q$ that do not occur in $F_\lambda$.
At rational masses, \eqref{eq:pointaccept} is an inequality between
rational numbers. The integer factors used for the mass domains in
Section~\ref{sec:cover} are part of the Zenodo deposit.

\section{Coefficient estimates and large masses}\label{sec:kernels}
The preceding lemma controls the angular parameters at fixed masses.
We now obtain bounds for variation in the masses and for the region
where the total mass is large. Both follow from coefficient estimates
for a single function with positive real part.

\subsection{Transfer of coefficient bounds}
Put
\[
 H_0(w)=\frac{1+w}{1-w},\qquad
 d_{k,\ell}(x)=(\partial_x-1)^\ell p_k(x).
\]
Thus $e^{-x}d_{k,\ell}(x)$ is the coefficient of $w^k$ in
$(-H_0(w))^\ell e^{-xH_0(w)}$.

\begin{lemma}[Coefficient transfer]\label{lem:transfer}
Let $x\geq0$ and let $\ell$ be a nonnegative integer. Suppose that
$H$ and $g$ are holomorphic on $\D$, with $H(0)=1$,
$\Rea H>0$, and $|g|\leq1$. If a polynomial $K$ agrees through
degree seven with $(-H_0)^\ell e^{-xH_0}$, then
\begin{equation}\label{eq:transfer}
 \left|[z^7](-H)^\ell e^{-xH}g\right|\leq\|K\|_{H^1}.
\end{equation}
In particular,
\begin{equation}\label{eq:h2}
 \left|[z^7](-H)^\ell e^{-xH}g\right|
 \leq B_\ell(x):=
 e^{-x}\left(\sum_{k=0}^7d_{k,\ell}(x)^2\right)^{1/2}.
\end{equation}
When $g=1$, the constant coefficient of $K$ may be chosen freely.
\end{lemma}
\begin{proof}
The map $\omega=(H-1)/(H+1)$ takes $\D$ into itself and vanishes
at zero, and $H=H_0\circ\omega$. The difference between $K$ and
the prescribed function is $O(w^8)$ at zero. After composition
with $\omega$ and multiplication by $g$, it is still $O(z^8)$.
Thus the coefficient in \eqref{eq:transfer} equals $[z^7](K\circ\omega)g$.

Let $u$ be the harmonic function on the disk with boundary values
$|K|$. Subharmonicity gives $|K|\leq u$, and $u\circ\omega$ is
a harmonic majorant of $|K\circ\omega|$. Its mean on each circle
centered at zero is $u(\omega(0))=\|K\|_1$. Consequently
$\|K\circ\omega\|_{H^1}\leq\|K\|_1$.
Multiplication by $g$ cannot increase this norm, and Cauchy's
coefficient formula bounds each Taylor coefficient by the $H^1$
norm. This proves \eqref{eq:transfer}.

Choose $K=e^{-x}\sum_{k=0}^7d_{k,\ell}(x)w^k$ and use
$\|K\|_1\leq\|K\|_2$ to obtain \eqref{eq:h2}.
If $g=1$, a change in the constant term remains constant after
composition and hence does not change the seventh coefficient.
\end{proof}

The functions $B_\ell$ are nonincreasing. The following identity
proves this without locating the critical points of a polynomial.
\begin{lemma}[Monotonicity of coefficient norms]\label{lem:monotonicity}
Let $H_0(w)=(1+w)/(1-w)$, let $G$ be analytic in the disk and independent
of $x$, and put
\[
 c_k(x)=[w^k]\bigl(e^{-xH_0(w)}G(w)\bigr),\qquad 0\leq k\leq n.
\]
For real $x$ and every integer $n\geq0$,
\begin{equation}\label{eq:jet-dissipation}
 \frac{d}{dx}\sum_{k=0}^n|c_k(x)|^2
   =-2\left|\sum_{k=0}^n c_k(x)\right|^2\leq0.
\end{equation}
In particular, for every integer $\ell\geq0$ and
$d_{k,\ell}=(\partial_x-1)^\ell p_k$,
\[
 \frac{d}{dx}\left(e^{-2x}\sum_{k=0}^n d_{k,\ell}(x)^2\right)
   =-2e^{-2x}\left(\sum_{k=0}^n d_{k,\ell}(x)\right)^2.
\]
Hence each function $B_\ell$ in \eqref{eq:h2} is
nonincreasing on $[0,\infty)$.
\end{lemma}
\begin{proof}
Since $H_0(w)=1+2\sum_{j\geq1}w^j$, differentiation of a fixed
Taylor coefficient gives
\[
 c_k'=-c_k-2\sum_{j<k}c_j.
\]
Consequently
\[
 2\Rea\sum_{k=0}^n c_k'\overline{c_k}
 =-2\sum_{k=0}^n|c_k|^2
   -4\Rea\sum_{0\leq j<k\leq n}c_j\overline{c_k}
 =-2\left|\sum_{k=0}^nc_k\right|^2.
\]
Taking $G=(-H_0)^\ell$ gives $c_k=e^{-x}d_{k,\ell}$ and proves
the second identity.
\end{proof}

The sum includes the constant coefficient. An angular derivative
introduces weights on the coefficients, so its norm requires a
separate estimate.

\subsection{Polynomial completions}
The estimate $\|K\|_1\leq\|K\|_2$ can be improved by choosing
coefficients of $K$ above degree seven. A useful construction starts
with real matrices $L,M$ having eight rows. Define
\[
 P(z)=\sum_r\left(\sum_{i=0}^7L_{ir}z^i\right)^2
      -\sum_r\left(\sum_{i=0}^7M_{ir}z^i\right)^2,
 \qquad E=\|L\|_F^2+\|M\|_F^2,
\]
where $\|L\|_F^2=\sum_{i,r}|L_{ir}|^2$.
For prescribed coefficients $c_0,\ldots,c_7$, put
\[
 K(z)=P(z)+\sum_{k=0}^7(c_k-[z^k]P)z^k.
\]
Parseval's identity gives $\|A^2\|_1=\|A\|_2^2$ for every
analytic polynomial $A$. The triangle inequality therefore implies
\begin{equation}\label{eq:smallkernel}
 \|K\|_{H^1}\leq E+
       \left(\sum_{k=0}^7|c_k-[z^k]P|^2\right)^{1/2}.
\end{equation}
This estimates the $H^1$ norm using only polynomial coefficients.

For second mass derivatives, take $c_k=e^{-x}d_{k,2}(x)$.
On a closed interval centered at $x_0$, the variation of the vector
$c=(c_0,\ldots,c_7)$ is bounded by its derivative, whose norm is
$B_3(x)$. More generally, for an interval of radius $r$,
\[
 \sup_{|x-x_0|\leq r}\|c'(x)\|_2
 \leq\|c'(x_0)\|_2+
          r\sup_{|x-x_0|\leq r}\|c''(x)\|_2.
\]
Together with the monotonicity identity, these estimates extend
\eqref{eq:smallkernel} from the midpoint to the whole interval.

The finite estimates used below are obtained on 2,700 closed
intervals of width $1/700$ covering $[1/7,4]$.
For $x\geq4$, \eqref{eq:h2} and monotonicity give
$B_2(x)\leq B_2(4)<149147/200000$.
For fixed $H,g$ as in Lemma~\ref{lem:transfer}, write
$A(x)=[z^7]e^{-xH}g$. The resulting bounds for $|A''(x)|$ are
summarized by
\begingroup
\setlength{\arraycolsep}{4pt}
\[
\begin{array}{c|rrrrrr}
 x\geq &1/7&2/7&3/7&4/7&5/7&6/7\\\hline
 |A''(x)|\leq&9.875458&4.725215&3.875485&3.703567&3.014671&2.376301
\end{array}
\]
\endgroup
The proof of the mass-domain estimate uses the sharper bounds on
individual intervals, taking the maximum over intervals that meet
the relevant range of $x$.

To bound $a_7$ in terms of the total mass alone, take
$c_k=e^{-t}p_k(t)$ for $1\leq k\leq7$ and choose the constant
freely, as permitted by Lemma~\ref{lem:transfer}. This gives bounds
on 1,815 closed intervals covering $[6/7,69/20]$.
Each bound is at most $797919/10^6$; again, the individual values
are used in the interpolation argument.
The finite inequalities asserted in this subsection are included
in Proposition~\ref{prop:finite}, with their rational factors
specified in the Zenodo deposit.

\subsection{A uniform estimate for large total mass}
\begin{lemma}[Quadrature on a mass interval]\label{lem:quadrature}
Let $K_x$ be a polynomial whose coefficients depend continuously
differentiably on a real parameter $x\in[x_0-h/2,x_0+h/2]$.
For a positive integer $m$, put
\[
 Q_m(K)=\frac1m\sum_{j=0}^{m-1}|K(e^{2\pi i j/m})|.
\]
Then
\[
 \|K_x\|_{H^1}
 \le Q_m(K_{x_0})
 +\frac{\pi}{m}\|\partial_\theta K_{x_0}\|_{H^2}
 +\frac h2\sup_{|y-x_0|\le h/2}\|\partial_yK_y\|_{H^2}.
\]
\end{lemma}
\begin{proof}
Partition the circle into the $m$ arcs centered at the quadrature nodes.
On each arc, the fundamental theorem of calculus bounds the difference
between $|K_{x_0}|$ and its central value by the integral of
$|\partial_\theta K_{x_0}|$ between those points. Integrating and reversing
the order of integration gives an error at most
$\pi\|\partial_\theta K_{x_0}\|_{H^1}/m$.
This argument remains valid at polynomial zeros, since the modulus is
absolutely continuous. Finally,
\[
 \|K_x-K_{x_0}\|_{H^1}
 \le |x-x_0|\sup_y\|\partial_yK_y\|_{H^1},
\]
and $\|\cdot\|_{H^1}\le\|\cdot\|_{H^2}$ gives the result.
\end{proof}

Lemma~\ref{lem:monotonicity} bounds the mass derivatives in
Lemma~\ref{lem:quadrature} by endpoint values. For a polynomial of
degree seven, Parseval's identity gives
$\|\partial_\theta K\|_2\leq7\|K\|_2$.
For a polynomial of degree fifteen, the angular norm at the midpoint
is evaluated from its coefficients. The finite values used below
are rigorous upper bounds, with both quadrature errors included.

\begin{lemma}[Large total mass]\label{lem:tail}
For every positive Herglotz measure of mass $t\geq69/20$,
\begin{equation}\label{eq:tail345}
 |a_7|<14697/20000=0.73485<\beta.
\end{equation}
\end{lemma}
\begin{proof}
For $t\geq703/100$, subordination gives $|a_7|\leq B_0(t)$.
Lemma~\ref{lem:monotonicity} and the endpoint inequality
$B_0(703/100)<14697/20000$ prove the stated bound on that half-line.
Here $B_0$ includes the constant coefficient, as required for monotonicity.

On $[18/5,703/100]$ use
$R_t=e^{-t}\sum_{k=1}^7p_k(t)w^k$. Partition the interval
into 1,715 closed intervals of width $1/500$. At each midpoint,
the 4,096-point mean is below $0.72974$. The endpoint bounds
$B_0(18/5)<0.813$ and $B_1(18/5)<0.71$, together with monotonicity,
give $\|\partial_\theta R_t\|_2\leq7B_0(t)<5.691$ and
$\|\partial_t R_t\|_2\leq B_1(t)<0.71$. Thus
\[
 \|R_t\|_1<0.72974+5.691\pi/4096+0.71/1000
 <0.734814943<0.73485.
\]
On $[69/20,18/5]$, 150 closed intervals of width $1/1000$ use polynomials
with the prescribed coefficients one through seven and rational
free coefficients in degrees zero and eight through fifteen.
The free coefficients stay fixed on each interval. The 32,768-point
mean, the angular norm of the midpoint polynomial, and the uniform
mass derivative bound $B_1(69/20)$ bound their $H^1$ norms by
$0.7343$. These polynomials preserve the seventh coefficient under
composition. Combining the three intervals proves
\eqref{eq:tail345}. The finite comparisons in this proof are part of
Proposition~\ref{prop:finite}; the finite data are described in
Appendix~\ref{app:verification} and specified in the Zenodo deposit.
\end{proof}

\section{The bounded mass domains}\label{sec:cover}
It remains to treat $N=3,4,5,6$ with $t\leq69/20$.
After permuting the masses and their support parameters together,
assume $\lambda_1\leq\cdots\leq\lambda_N$. Define nonnegative
weighted gaps $u_1,\ldots,u_N$ by
\begin{equation}\label{eq:gaps}
 \lambda_j=\frac17+\sum_{k=1}^j\frac{u_k}{N-k+1},
 \qquad
 t=\frac N7+\sum_{j=1}^Nu_j.
\end{equation}
This parametrizes all ordered masses with $\lambda_j\geq1/7$,
including equal adjacent masses. The domain to be bounded is
\begin{equation}\label{eq:mass-slab}
 \frac7N-\frac N7\leq\sum_{j=1}^Nu_j
     \leq\frac{69}{20}-\frac N7,
 \qquad u_j\geq0.
\end{equation}
It is contained in the cube
\begin{equation}\label{eq:rootcube}
 C_N=[0,18/5-N/7]^N.
\end{equation}
We cover this cube by finitely many closed boxes and apply
interpolation on each box that meets \eqref{eq:mass-slab}.

\subsection{Second derivatives in a gap coordinate}
Fix $j$, put $m=N-j+1$, and write
\[
 H_k(z)=\frac{1+\alpha_kz}{1-\alpha_kz},\qquad
 H=\frac1m\sum_{k=j}^NH_k.
\]
The factorization
\begin{align*}
 f&=e^{-xH}g,\qquad x=m\lambda_j,\\
 g&=\exp\!\left[-\sum_{k<j}\lambda_kH_k
              -\sum_{k\geq j}(\lambda_k-\lambda_j)H_k\right]
\end{align*}
has $|g|\leq1$, because all the masses in its exponent are
nonnegative. Varying $u_j$ while fixing the other gaps changes $x$
at rate one and leaves $H$ and $g$ fixed. Thus
\[
 \partial_{u_j}^2a_7=[z^7]H^2e^{-xH}g.
\]
Lemma~\ref{lem:transfer} therefore bounds this derivative on a box
using the interval of values of $x=m\lambda_j$ between its two
opposite corners. The estimate is uniform in the support parameters.

\subsection{Interpolation on a box}
Write the box as $u=a+w\odot s$, where $s\in[0,1]^N$ and
$w\odot s=(w_1s_1,\ldots,w_Ns_N)$.
For values $b_\varepsilon$ at its vertices, define their multilinear
interpolant by
\[
 \operatorname{ML}(b)(s)=
 \sum_{\varepsilon\in\{0,1\}^N}b_\varepsilon
       \prod_{j=1}^Ns_j^{\varepsilon_j}(1-s_j)^{1-\varepsilon_j}.
\]
The factors in this sum are nonnegative and have sum one.

\begin{lemma}[Interpolation estimate]\label{lem:interpolation}
Suppose that $|\partial_{u_j}^2a_7|\leq M_j$ throughout the box,
and that $|a_7(a+w\odot\varepsilon)|\leq b_\varepsilon$ at every
vertex. Then
\begin{equation}\label{eq:multilinear}
 |a_7(a+w\odot s)|\leq q(s):=
 \operatorname{ML}(b)(s)+\sum_{j=1}^Nc_js_j(1-s_j),
 \qquad c_j=\frac{M_jw_j^2}{2}.
\end{equation}
\end{lemma}
\begin{proof}
In one variable, the integral remainder for linear interpolation
bounds the error by $M_jw_j^2s_j(1-s_j)/2$. This bound also holds
for a complex-valued function, by taking absolute values in the
integral remainder. Apply the interpolation operators in the $N$
coordinates successively. Their nonnegative weights preserve each
error bound, and the resulting errors sum to the second term in
\eqref{eq:multilinear}. The absolute value of the final interpolant
is at most $\operatorname{ML}(b)$.
\end{proof}

The support parameters are fixed during this argument, but both
the vertex bounds and the derivative bounds hold uniformly in
them. No choice of maximizing parameters is needed.

To bound $q$ on the unit cube, use the elementary identity
\[
\begin{aligned}
 T(a,b,c)&=\max_{0\leq s\leq1}\{(1-s)a+sb+cs(1-s)\}\\
 &=\begin{cases}
 \max(a,b),&c=0\text{ or }|a-b|\geq c,\\
 (a+b)/2+c/4+(a-b)^2/(4c),&|a-b|<c,
 \end{cases}
\end{aligned}
\]
valid for real $a,b$ and $c\geq0$.
For fixed $c$, this function is convex and nondecreasing in $a,b$,
and satisfies $T(a+d,b+d,c)=T(a,b,c)+d$.
The endpoint values in one coordinate are convex combinations of
pairs of vertex values. Jensen's inequality therefore bounds the
maximum in that coordinate by replacing each pair with its $T$
value. Repeating this operation gives a rational upper bound for
$q$. Different orders of elimination may be used, and the least
of the resulting upper bounds remains valid.

The slab restriction \eqref{eq:mass-slab} can sharpen this bound.
Set
\[
 t_0=\frac N7+\sum_ja_j,\qquad
 R_- =\frac7N-t_0,\qquad R_+=\frac{69}{20}-t_0.
\]
For a rational number $\rho$ and a point of the box in the slab,
\begin{equation}\label{eq:slab}
 q(s)\leq\rho R+
    \max_{[0,1]^N}\{q(s)-\rho w\cdot s\},
 \qquad
 R=\begin{cases}R_+,&\rho\geq0,\\R_-,&\rho<0.\end{cases}
\end{equation}
The maximum on the right is bounded by the same argument, with
vertex values $b_\varepsilon-\rho w\cdot\varepsilon$.
These values may be negative; convexity and monotonicity of $T$
still apply. Every choice of $\rho$ gives a valid inequality.

At each vertex we use the smaller of the applicable bounds from
Sections~\ref{sec:phase} and \ref{sec:kernels}. When neither bound
is available, Cauchy's estimate gives the universal bound $1$.
Thus all vertices have valid bounds, including vertices outside
the slab.

\subsection{A finite cover}
Starting with $C_N$, bisect one coordinate of a box at its rational
midpoint, retaining both closed halves. The deposited data
specify a finite sequence of such bisections. Each final box has
one of the following properties:
\begin{enumerate}[label=\textup{(\roman*)}]
\item its largest total mass is less than $7/N$;
\item its smallest total mass is at least $69/20$;
\item its upper bound from \eqref{eq:slab} is less than $\beta$.
\end{enumerate}
A box of the first kind does not meet the admissible domain.
Lemma~\ref{lem:tail} applies to the second kind, and
Lemma~\ref{lem:interpolation} applies to the third.
Each bisection preserves the union of its two children, so induction
proves that the final boxes cover $C_N$, including all their
boundary faces.

\begin{proposition}[Finite estimates]\label{prop:finite}
The finite data described in Appendix~\ref{app:verification}
give the following bounds and subdivisions.
\begin{enumerate}[label=\textup{(\roman*)}]
\item At $66{,}631$ rational mass vectors, the factors satisfy
\eqref{eq:pointaccept} with their specified rational coefficient bounds.
\item All $2{,}700$ second-derivative cells and $1{,}815$ scalar cells
satisfy \eqref{eq:smallkernel} with their specified rational bounds.
The endpoint comparisons and all $1{,}865$ finite quadrature cells
used in Lemma~\ref{lem:tail} satisfy its stated inequalities.
\item For each $N=3,4,5,6$, the subdivision partitions the cube
in \eqref{eq:rootcube} by closed rational midpoint bisections. Every
final box has one of the three properties above, with every
vertex and derivative bound linked to the data in \textup{(i)} and
\textup{(ii)} or to the universal bound $1$.
\end{enumerate}
The sizes of these subdivisions are given in
Table~\ref{tab:certificates}.
\end{proposition}
\begin{proof}[Finite verification]
The first assertion consists of rational comparisons
\eqref{eq:pointaccept}, with all Fourier coefficients included.
For the second, \eqref{eq:smallkernel} is bounded on each closed
interval using the coefficient variation estimates of
Section~\ref{sec:kernels}. The quadrature bounds include both
error terms in Lemma~\ref{lem:quadrature}.
For the third, the rational bisections determine every box exactly;
its vertex and derivative bounds give the corresponding inequality
\eqref{eq:slab}. All final boxes satisfy one of the three stated
conditions. Appendix~\ref{app:verification} describes the arithmetic
for these finite comparisons. Their full factors, bounds, and
subdivisions are included in the Zenodo deposit.
\end{proof}

\begin{table}[ht]
\centering\small
\begin{tabular}{@{}rrrrrrr@{}}
\toprule
$N$ & $\dim v$ & Points & Splits & Interpolation & Tail & Empty\\
\midrule
3&25&257&215&186&29&1\\
4&55&1,926&1,041&915&127&0\\
5&105&10,262&3,819&3,489&331&0\\
6&182&54,186&13,022&12,252&771&0\\
\bottomrule
\end{tabular}
\caption{Finite subdivisions of the mass domains. The Points column
counts rational mass vectors; the last three columns count final boxes
bounded by interpolation, by the large-mass estimate, or excluded
by the lower mass bound. The subdivisions contain $36{,}198$ boxes
in total, including those that are further bisected.}
\label{tab:certificates}
\end{table}
The largest interpolation bound among the final boxes is approximately
$0.734998509736$. Its comparison with $\beta=0.735$ is made as an
exact rational inequality.

\begin{theorem}[Deficient-support exclusion]\label{thm:exclusion}
For $1\leq N\leq6$, all functions \eqref{eq:atoms} with
$\lambda_j\geq1/7$ and $t\geq7/N$ satisfy $|a_7|<\beta$.
\end{theorem}
\begin{proof}
For $N=1,2$, the total-mass bound gives $t\geq7/2>69/20$ and
Lemma~\ref{lem:tail} applies. For $N=3,4,5,6$,
Proposition~\ref{prop:finite} gives the bound on the finite slab,
and Lemma~\ref{lem:tail} gives the bound for larger masses.
Simultaneous permutation of masses and support parameters removes
the ordering convention.
\end{proof}

\section{Proof of the main theorem}\label{sec:conclusion}
We now combine the deficient-support exclusion with the entropy bound.

\begin{proof}[Proof of Theorem~\ref{thm:fullsupport}]
A seventh-coefficient extremizer exists and has at most seven
distinct positive atoms by Section~\ref{sec:mass}.
If it had $N\leq6$ atoms, Lemmas~\ref{lem:individual}
and~\ref{lem:total} would give $\lambda_j\geq1/7$ and $t\geq7/N$.
Theorem~\ref{thm:exclusion} would then imply
$M_7=|a_7|<147/200<2/e$, contrary to \eqref{eq:model}.
Hence every extremizer has exactly seven atoms.
\end{proof}

\begin{proof}[Proof of Theorem~\ref{thm:seventh}]
Apply Corollary~\ref{cor:fullsupport} with $n=7$, using
Theorem~\ref{thm:fullsupport}. The Taylor expansion in
\eqref{eq:model} verifies equality for every stated rotation.
\end{proof}

\appendix
\section{The finite inequalities}\label{app:verification}
Proposition~\ref{prop:finite} is the finite part of the proof.
Its inputs are rational masses, integer matrix factors, rational
interval endpoints, and finite descriptions of the subdivisions.
The analytic arguments in Sections~\ref{sec:phase}--\ref{sec:cover}
establish why the corresponding inequalities suffice.
The complete certificate data, verification programs, manifests,
and verification records form a separate Zenodo deposit.
This appendix describes the data, the arithmetic, and the
verification procedure.

\subsection{Data inventory}
At each of the $66{,}631$ rational mass vectors in
Table~\ref{tab:certificates}, the data specify the integer matrix
factor, its rational scale, and the coefficient bound used in
\eqref{eq:pointaccept}. All Fourier coefficients are included in
the residual comparison. The four mass subdivisions record the
rational midpoint bisections, the vertex and derivative bounds,
and the slab multipliers. Together they have $36{,}198$ nodes.
A binary tree with $s$ splits has $2s+1$ nodes and $s+1$ leaves;
these identities agree with the counts in the table.

The interval data comprise:
\begin{itemize}
\item $2{,}700$ second-derivative cells of width $1/700$
covering $[1/7,4]$, with their rational polynomial completions
and bounds;
\item $1{,}815$ scalar cells covering $[6/7,69/20]$, with
their individual rational bounds;
\item $150$ quadrature cells covering $[69/20,18/5]$, each
using $32{,}768$ angular nodes;
\item $1{,}715$ quadrature cells covering $[18/5,703/100]$,
each using $4{,}096$ angular nodes.
\end{itemize}
Thus the interval comparisons use $4{,}515$ derivative and
scalar cells and $1{,}865$ quadrature cells. The remaining
large-mass range follows from the endpoint and monotonicity
estimates in Lemma~\ref{lem:tail}.

\subsection{Rational point bounds}
At a rational mass vector, every coefficient in
\eqref{eq:phasepoly} is rational. To see explicitly how denominators
may be cleared, write $\lambda_i=N_i/D$. Then
\[
 P_{ik}=D^k k!\,p_k(N_i/D)
\]
is an integer by \eqref{eq:pk}. For $|k|=7$, the numerator of
$\prod_i p_{k_i}(\lambda_i)$ over the denominator $D^7 7!$ is
\[
 \binom{7}{k_1,\ldots,k_N}\prod_iP_{i,k_i}.
\]
The Laurent polynomial $Q$ also has rational coefficients, since
it is given by an integer factor and a rational scale.
Thus every term in $U$ and $bE_{128}(t)$ can be compared exactly.
No estimate of the eigenvalues of a matrix is involved: the
nonnegativity of $Q$ follows directly from its sum-of-squares form.

\subsection{Bounds on intervals}
For the completions in \eqref{eq:smallkernel}, the matrices are
rational, so their squared norms and polynomial coefficients are
rational as well. The remaining expressions involve $e^{-x}$
and square roots of nonnegative quantities. Taylor's integral
remainder gives, for $x\geq0$,
\[
 \sum_{k=0}^{99}\frac{(-x)^k}{k!}
 \leq e^{-x}\leq
 \sum_{k=0}^{100}\frac{(-x)^k}{k!}.
\]
At rational endpoints these are rational enclosures. An upper bound
$u$ for $\sqrt r$, with $u\geq0$ and $r\geq0$, is verified by
$u^2\geq r$. The variation estimates in
Section~\ref{sec:kernels} then give bounds throughout each closed
interval, including its endpoints.

The finite quadrature means in Lemma~\ref{lem:tail} use rigorous
real and complex interval enclosures. The angular and mass-variation
errors in Lemma~\ref{lem:quadrature} are added to the upper endpoint
of each enclosure. The resulting upper bounds are compared with
the rational constants stated in Lemma~\ref{lem:tail}.

\subsection{Coverage of the mass domains}
Each subdivision starts from the exact cube \eqref{eq:rootcube}.
A bisection records a coordinate and its rational midpoint;
the two resulting closed boxes have union equal to their parent.
The bounds assigned to a final box are obtained from its vertices,
its intervals of grouped masses, and the chosen rational multiplier
$\rho$ in \eqref{eq:slab}.
The three final-box conditions in Section~\ref{sec:cover} are
checked as strict rational inequalities, except for the weak
comparison placing a box in the large-mass region, where
Lemma~\ref{lem:tail} already gives a strict coefficient bound.
This verifies both the local estimates and the coverage of the
entire admissible set, including the common faces of adjacent boxes.

\subsection{Verification and data availability}
The Zenodo distribution contains the unchanged certificate
archive \path{certificate-data.tar.gz}, the portable sources in
\path{verifier-source.tar.gz}, and the successful run records in
\path{verification-records.tar.gz}. It also includes the compressed
original data manifest, a file inventory, and
\path{DATA_AVAILABILITY.json}, which records the archive's exact
size and SHA-256 digest.

Extract the three archives into a common working directory.
The combined driver is
\path{portable_verification/CHECK_REDUCED.py}. Its dependencies
are Python with \texttt{python-flint}, Node.js, a C++17 compiler
with \texttt{\_\_int128} support, and zlib. Use paths without
spaces or URI-escaped characters, as required by this release.
With assertions enabled, run
\begin{quote}\small\ttfamily
python3 -B portable\_verification/CHECK\_REDUCED.py \textbackslash\\
\hspace*{1em}--archive seventh\_coefficient\_certificate\_20261029 \textbackslash\\
\hspace*{1em}--output fresh\_replay
\end{quote}
The output directory must be new and outside the certificate
directory. The driver reads the frozen data without modification
and writes its commands, logs, source snapshots, and receipts
to the output directory.

The accompanying \path{verification/REDUCED_CHAIN.json} records
a complete successful run: $66{,}709$ manifest-listed archive
files, $66{,}631$ point factors, $4{,}515$ interval cells,
$36{,}198$ subdivision nodes, and $1{,}865$ quadrature cells,
with no unresolved leaves. The source and data hashes identify
the exact inputs to those checks. File hashes establish identity;
the arithmetic and coverage checks establish the finite estimates.

% DATA_AVAILABILITY_BEGIN
The certificate data, verification sources, and run records are
associated with the Zenodo identifier
\url{https://doi.org/10.5281/zenodo.23816750}.
% DATA_AVAILABILITY_END

\section{The simple-root entropy inequality}\label{app:entropy}
We give the simple-root case of the entropy inequality of Lei and
Zhang \cite[Section~6]{LeiZhang}, including their quantitative
remainder. This is the form used in Section~\ref{sec:entropy}.

\begin{lemma}[Lei--Zhang, simple-root case]\label{lem:entropy}
Let $p$ have degree $n\geq2$ and distinct zeros on $\T$.
After a unimodular normalization write
$p=p^*=\sum_{j=0}^n c_jz^j$, where reflection is relative to degree $n$.
With $A=\|p\|_2^2$ and
$\Gamma=\sum_{j=1}^{n-1}j(n-j)|c_j|^2/n^2$, one has
\[
 \int_{\T}|p|^2\log|p|^2\,dm
 \geq A\bigl(1+\log(A/2)\bigr)+\frac{2\Gamma}{n(n-1)}.
\]
Equality in the bound without the last term requires
$p=c(\omega+z^n)$, $|\omega|=1$; these polynomials attain it.
\end{lemma}
\begin{proof}
Put $q=p-zp'/n$, with $q^*$ reflected relative to degree $n$
although $\deg q<n$. Self-inversivity gives $p=q+q^*$ and
$q^*=zp'/n$. Inside the disk,
$\Rea(zp'/p)=\sum_{p(\tau)=0}\Rea(z/(z-\tau))<n/2$,
so $q$ cannot vanish there. A boundary zero of $q$ would also
annihilate $q^*$ and then both $p$ and $p'$, contradicting simplicity.
Since $q$ has no zeros in the closed disk, $r=q^*/q$ is holomorphic on a neighborhood of the
closed disk, has modulus one on $\T$, and satisfies $r(0)=0$.
In particular, $|r|\leq1$ in $\D$.

The measure $d\nu=|1+r|^2dm/2$ has total mass one, since
$|r|=1$ on $\T$ and $\int_\T r\,dm=r(0)=0$.
Moreover, $p=q(1+r)$ gives $\int|q|^2\,d\nu=A/2$.
Jensen's inequality for the convex function $x\log x$ therefore gives
\[
 J:=\int |p|^2\log|q|^2\,dm
   =2\int|q|^2\log|q|^2\,d\nu
   \geq A\log(A/2).
\]
We next estimate the part of the entropy involving $1+r$.
For a holomorphic function $g(z)=\sum_{j\geq0}g_jz^j$ with
$g(0)=0$, define
$S(g)=\sum_{j=1}^{n-1}(n-j)|g_j|^2/j$.
If $\varphi$ is analytic with $|\varphi|\leq1$ in
the disk, let $P_\ell$ denote Taylor truncation through degree
$\ell$. Projection and multiplication by $\varphi$ are contractions
in $H^2$, so
\[
 \norm{P_\ell(\varphi g)}_2
 =\norm{P_\ell(\varphi P_\ell g)}_2
 \leq\norm{P_\ell g}_2.
\]
Put $w_j=(n-j)/j$ for $1\leq j<n$ and $w_n=0$. Since
\[
 S(g)=\sum_{\ell=1}^{n-1}(w_\ell-w_{\ell+1})\|P_\ell g\|_2^2,
\]
summing the squared inequalities with these positive differences proves
$S(\varphi g)\leq S(g)$.

Use the inner product $\langle u,v\rangle=\int u\bar v\,dm$
and let $M_k=\langle r^kq,q\rangle$.
The formulas $q_j=(n-j)c_j/n$ and $(q^*)_j=jc_j/n$ imply
$M_1=S(q^*)=\Gamma$.
The dual weighted norm also satisfies
$\sum_{j=1}^{n-1}|q_j|^2/w_j=\Gamma$.
For $g=r^{k-1}q^*=r^kq$, the preceding contraction gives
$S(g)\leq\Gamma$. Since $g(0)=0$ and $q_n=0$, weighted
Cauchy--Schwarz gives
\[
 |M_k|^2=\left|\sum_{j=1}^{n-1}g_j\overline{q_j}\right|^2
 \leq S(g)\sum_{j=1}^{n-1}\frac{|q_j|^2}{w_j}\leq\Gamma^2.
\]
Also $M_k=0$ for $k\geq n$, since $r^kq$ vanishes to order at least
$k$ and $\deg q<n$.

On $|w|=1$, the absolutely convergent Fourier identity is
\[
 |1+w|^2\log|1+w|^2
 =2+\tfrac32(w+\bar w)
 +\sum_{k=2}^{\infty}\frac{2(-1)^k(w^k+\bar w^k)}{k(k^2-1)}.
\]
It follows by multiplying the Abel series for $\log|1+w|^2$ by
$2+w+\bar w$; the resulting summable coefficients justify the
boundary identity and substitution $w=r$.
Since $A=2\|q\|_2^2+2\Gamma$, writing
$I=\int |p|^2\log|p/q|^2\,dm$ gives
\[
 I-A=\Gamma+4\sum_{k=2}^{n-1}
           \frac{(-1)^k\Rea M_k}{k(k^2-1)}
 \geq\frac{2\Gamma}{n(n-1)},
\]
because $\sum_{k=2}^{n-1}1/[k(k^2-1)]=1/4-1/[2n(n-1)]$.
Since the entropy is $I+J$, adding the bound for $J$ proves the
claim. Equality in the unstrengthened inequality
forces $\Gamma=0$, hence all intermediate coefficients vanish;
self-inversivity gives the stated form. Its integral follows from
the constant Fourier coefficient $2$ in the displayed identity.
\end{proof}

\section{Use of AI tools}\label{app:ai}
ChatGPT (OpenAI GPT-6.1 Sol) was used to assist with code development,
checking arguments, and editing the manuscript. The author is
responsible for all mathematical content and code developed for
this work.

\end{document}